\newcount\ite\ite=1\def\0{\global\ite=1\1}
\def\1{\item{\rm(\romannumeral\the\ite)}\advance\ite1\quad}
\def\phi{\varphi}

\documentclass[12pt,twoside]{article}
\usepackage{amssymb}

\font\teneufm=eufm10 scaled \magstep1
\font\seveneufm=eufm7 scaled \magstep1
\font\fiveeufm=eufm5  scaled \magstep1
\newfam\eufmfam

\textfont\eufmfam=\teneufm
\scriptfont\eufmfam=\seveneufm
\scriptscriptfont\eufmfam=\fiveeufm
\def\frak#1{{\fam\eufmfam\relax#1}}

\newfam\msbfam
\font\tenmsb=msbm10 scaled \magstep1  \textfont\msbfam=\tenmsb
\font\sevenmsb=msbm7 scaled \magstep1 \scriptfont\msbfam=\sevenmsb
\font\fivemsb=msbm5 scaled \magstep1  \scriptscriptfont\msbfam=\fivemsb

\makeatletter
\def\blfootnote{\xdef\@thefnmark{}\@footnotetext}
\makeatother

\def\dd#1{\raise1.5pt\hbox{$\,\partial\!$}/\raise-2.5pt\hbox{$\!\partial#1\,$}}

\def\tilde{\widetilde}

\def\5#1{{\mathcal #1}}

\def\CC{{\mathbb C}}

\def\PP{{\mathbb P}}

\def\FF{{\mathbb F}}

\def\ra{\rightarrow}

\def\GL{\mathop{\rm GL}\nolimits}

\def\Ann{\mathop{\rm Ann}\nolimits}

\def\emb{\mathop{\rm emb}\nolimits}
 \def\HollowBoxx #1#2#3{{\dimen0=#1 \advance\dimen0 by -#2
       \dimen1=#1 \advance\dimen1 by #3
        \vrule height 0pt depth #3 width #2
       \hskip -#3
       \vrule height #1 depth #3 width #3}}
 \def\LeftContraction{\mathord{\kern1.45pt \HollowBoxx{6pt}{3.5pt}{.4pt}}\,}

 \def\HollowBox #1#2#3{{\dimen0=#1 \advance\dimen0 by -#3
       \dimen1=#1 \advance\dimen1 by #3
        \vrule height #1 depth #3 width #3
        \vrule height 0pt depth #3 width #2
        \hskip -#3}}
 \def\RightContraction{\mathord{\, \HollowBox{6pt}{3.1pt}{.4pt}} \kern1.6pt}

\def\qed{{\hfill $\Box$}}
\newtheorem{theorem}{THEOREM}[section]
\newtheorem{corollary}[theorem]{Corollary}
\newtheorem{lemma}[theorem]{Lemma}

\newtheorem{remark}[theorem]{Remark}
\newtheorem{proposition}[theorem]{Proposition}

\begin{document}

\begin{center}
{\Large \bf On a New Criterion for Isomorphism\\
\vspace{0.3cm}
of Artinian Gorenstein Algebras}
\medskip\\
\normalsize A. V. Isaev
\end{center}

\begin{quotation} 
{\small \sl \noindent 
To every Gorenstein algebra $A$ of finite vector space dimension greater than 1 over a field $\FF$ of characteristic zero, and a linear projection $\pi$ on its maximal ideal 
${\mathfrak m}$ with range equal to the annihilator $\Ann({\mathfrak m})$ of ${\mathfrak m}$, one can associate a certain algebraic hypersurface $S_{\pi}\subset{\mathfrak m}$, which is the graph of a polynomial map $P_{\pi}:\ker\pi\ra\Ann({\mathfrak m})\simeq\FF$. Recently, in {\rm\cite{FIKK}}, {\rm\cite{FK}} the following surprising criterion was obtained: two Gorenstein algebras $A$, $\tilde A$ are isomorphic if and only if any two hypersurfaces $S_{\pi}$ and $S_{\tilde\pi}$ arising from $A$ and $\tilde A$, respectively, are affinely equivalent. The proof is indirect and relies on a CR-geometric argument. In the present paper we give a short algebraic proof of this statement. We also compare the polynomials $P_{\pi}$ with Macaulay's inverse systems. Namely, we show that the restrictions of $P_{\pi}$ to certain subspaces of $\ker\pi$ are inverse systems for $A$.}
\end{quotation}

\thispagestyle{empty}

\pagestyle{myheadings}
\markboth{A. V. Isaev}{Criterion for Isomorphism of Artinian Gorenstein Algebras}

\setcounter{section}{0}

\section{Introduction}\label{intro}
\setcounter{equation}{0}

We consider Artinian Gorenstein algebras over a field $\FF$ of characteristic zero. Recall that a local commutative associative algebra $A$ of finite vector space dimension greater than 1 is Gorenstein if and only if the annihilator $\Ann({\mathfrak m})$ of its maximal ideal ${\mathfrak m}$ (i.e., the socle of $A$) is 1-dimensional (see, e.g. \cite{Hu}). Gorenstein algebras frequently occur in various areas of mathematics and its applications to physics (see, e.g. \cite{B}, \cite{L}). In the case $\FF=\CC$, in our earlier paper \cite{FIKK} we found a surprising criterion for two such algebras to be isomorphic. The criterion was given in terms of a certain algebraic hypersurface $S_{\pi}\subset{\mathfrak m}$ associated to a linear projection $\pi$ on ${\mathfrak m}$ with range $\Ann({\mathfrak m})$. The hypersurface $S_{\pi}$ passes through the origin and is the graph of a polynomial map $P_{\pi}:\ker\pi\ra \Ann({\mathfrak m})\simeq\CC$. In \cite{FIKK} we showed that for $\FF=\CC$ two Gorenstein algebras $A$, $\tilde A$ are isomorphic if and only if any two hypersurfaces $S_{\pi}$ and $S_{\tilde\pi}$ arising from $A$ and $\tilde A$, respectively, are affinely equivalent, that is, there exists a bijective affine map ${\mathcal A}: {\mathfrak m}\ra\tilde {\mathfrak m}$ such that ${\mathcal A}(S_{\pi})=S_{\tilde\pi}$.

We stress that the above criterion for isomorphism of Gorenstein algebras differs from the well-known criterion in terms of Macaulay's inverse systems, where to an algebra $A$ one also associates a certain polynomial (see, e.g. Proposition 2.2 in article \cite{ER} and references therein). Indeed, as we show in Section \ref{invsys}, in general none of the polynomials $P_{\pi}$ is an inverse system for $A$. The advantage of the method of \cite{FIKK} lies in the fact that affine equivalence for the hypersurfaces $S_{\pi}$ and $S_{\tilde\pi}$ is much easier to understand than the type of equivalence required for inverse systems (see formula (7) in \cite{ER}). Indeed, in order to verify affine equivalence one can use methods of affine geometry, which makes this approach convenient in applications. In particular, it has proved to be quite useful for establishing isomorphism of isolated quasi-homogeneous hypersurface singularities (see e.g. \cite{EI}, \cite{I1}), since by the Mather-Yau theorem (see \cite{MY}) any such singularity is completely determined by its Milnor algebra, which is Gorenstein.

At the time of writing the article \cite{FIKK} we could not find an algebraic proof of our criterion, but, quite unexpectedly, discovered a proof based on the CR-geometry of certain real codimension two quadrics in complex space. The hypersurface $S_{\pi}$ can be introduced for a Gorenstein algebra over any field $\FF$ of characteristic zero, and therefore it is natural to attempt to extend the result of \cite{FIKK} to all such algebras. Indeed, an extension of this kind was obtained in recent paper \cite{FK}. The proof proceeds by reducing the case of an arbitrary field to the case $\FF=\CC$, and  therefore the basis of the method of \cite{FK} remains the geometric idea of \cite{FIKK}. It would be very desirable, however, to find a purely algebraic argument directly applicable to any field. In this paper we propose such an argument, in which neither CR-geometry nor reduction to the case $\FF=\CC$ is required. 
  
The paper is organized as follows. Section \ref{prelim} contains necessary preliminaries and the precise statement of the criterion in Theorem \ref{equivalence}. Our proof of Theorem \ref{equivalence} is given in Section \ref{proofs1}. In Section \ref{app} we demonstrate how powerful this result can be in applications. Namely, we apply it to a one-parameter family $A_t$ of 15-dimensional Gorenstein algebras. While directly finding all pairwise isomorphic algebras in the family $A_t$ seems to be quite hard, this problem is easily solved with the help of Theorem \ref{equivalence}. Finally, in Section \ref{invsys} we establish a connection between the polynomials $P_{\pi}$ and Macaulay's inverse systems for any Artinian Gorenstein algebra $A$. Namely, we show that if $P_{\pi}$ is regarded as a map from $\ker\pi$ to $\FF$ (in which case we call it a {\it nil-polynomial}), then the restriction of $P_{\pi}$ to any subspace of $\ker\pi$ that forms a complement to ${\mathfrak m}^2$ in ${\mathfrak m}$ is an inverse system for $A$ (see Theorem \ref{main}). All these subspaces have dimension equal to the embedding dimension $\emb\dim A$ of $A$, thus nil-polynomials can be viewed as certain extensions of inverse systems to spaces of dimension greater than $\emb\dim A$. The additional terms of nil-polynomials arising from such extensions have turned out to be rather useful in applications.

{\bf Acknowledgements.} We would like to thank Wilhelm Kaup and Nikolay Kruzhilin for useful discussions. This work is supported by the Australian Research Council.

\section{Preliminaries}\label{prelim}
\setcounter{equation}{0}
    
Let $A$ be a Gorenstein algebra of finite vector space dimension greater than 1 over a field $\FF$ of characteristic zero, and ${\mathfrak m}$ the maximal ideal of $A$. Let $\pi$ be a linear projection on $A$ with range $\Ann({\mathfrak m}):=\{u\in{\mathfrak m}: u\,{\mathfrak m}=0\}$ and kernel containing the identity element ${\bf 1}$ of $A$ (we call such projections {\it admissible}). Consider the following $\Ann({\mathfrak m})$-valued bilinear form on $A$:
\begin{equation}
b_{\pi}(a,c):=\pi(ac), \quad a,c\in A.\label{bpi}
\end{equation}
It is well-known that the form $b_{\pi}$ is non-degenerate (see, e.g. \cite{He}, p. 11). Further, let $\PP(A)$ be the projectivization of $A$ and
$$
Q_{\pi}:=\Bigl\{[a]\in\PP(A):b_{\pi}(a,a)=0\Bigr\},
$$
a projective quadric, where $[a]$ denotes the point of $\PP(A)$ represented by $a\in A$. The inclusion ${\mathfrak m}\subset A$ induces the inclusion $\PP({\mathfrak m})\subset \PP(A)$, and we think of $\PP({\mathfrak m})$ as the hyperplane at infinity in $\PP(A)$. We then identify ${\bf 1}+{\mathfrak m}\subset A$ with the affine part of $\PP(A)$ and introduce the corresponding affine quadric $Q'_{\pi}:=Q_{\pi}\cap ({\bf 1}+{\mathfrak m})$.

Let $\exp: {\mathfrak m}\ra {\bf 1}+{\mathfrak m}$ be the exponential map
$$
\displaystyle\exp(u):=\sum_{m=0}^{\infty}\frac{1}{m!}u^m,
$$
where $u^0:={\bf 1}$. By Nakayama's lemma the maximal ideal ${\mathfrak m}$ of $A$ is a nilpotent algebra, and therefore the above sum is in fact finite, with the highest-order term corresponding to $m=\nu$, where $\nu\ge 1$ is the nil-index of ${\mathfrak m}$ (i.e. the largest integer $\mu$ for which ${\frak m}^{\mu}\ne 0$). Thus, the exponential map is a polynomial transformation. It is bijective with the inverse given by the polynomial transformation
$$
\log({\bf 1}+u):=\sum_{m=1}^{\nu}\frac{(-1)^{m+1}}{m}u^m,\quad u\in{\mathfrak m}.
$$

We now define
$$
S_{\pi}:=2\log(Q_{\pi}')\subset{\mathfrak m}. 
$$
Letting ${\mathcal K}:=\ker\pi\cap{\mathfrak m}$, it is easily seen that $S_{\pi}$ is the graph of the polynomial map $P_{\pi}:{\mathcal K}\ra\Ann({\mathfrak m})$ of degree $\nu$ given by
\begin{equation}
P_{\pi}(u):=-\pi(\exp(u))=-\pi\left(\sum_{m=2}^{\nu}\frac{1}{m!}u^m\right),\quad u\in{\mathcal K}.\label{poly}
\end{equation}
It then follows that $S_{\pi}$ is an algebraic hypersurface in ${\mathfrak m}$ passing through the origin. Numerous examples of such hypersurfaces explicitly computed for particular algebras can be found in \cite{FIKK}, \cite{FK}, \cite{I1}, \cite{EI} (see also\linebreak Section \ref{app} below). 

We will now state the criterion for isomorphism of Gorenstein algebras obtained in \cite{FIKK}, \cite{FK}.

\begin{theorem}\label{equivalence}\sl Let $A$, $\tilde A$ be Gorenstein algebras of finite vector space dimension greater than 1, and $\pi$, $\tilde\pi$ admissible projections on $A$, $\tilde A$, respectively. Then $A$, $\tilde A$ are isomorphic if and only if the hypersurfaces $S_{\pi}$, $S_{\tilde\pi}$ are affinely equivalent. Moreover, if ${\mathcal A}:{\mathfrak m}\ra\tilde{\mathfrak m}$ is an affine equivalence such that ${\mathcal A}(S_{\pi})=S_{\tilde\pi}$, then the linear part of ${\mathcal A}$ is an algebra isomorphism between ${\mathfrak m}$ and $\tilde{\mathfrak m}$.
\end{theorem}

\section{Proof of Theorem \ref{equivalence}}\label{proofs1}
\setcounter{equation}{0}

The necessity implication follows from Proposition 2.2 of \cite{FIKK} (which works over any field of zero characteristic). Indeed, the proposition states that if $A$ is a Gorenstein algebra of dimension greater that 1 and $\pi_1$, $\pi_2$ are admissible projections on $A$, then $S_{\pi_{{}_1}}=S_{\pi_{{}_2}}+x$  for some $x\in{\mathfrak m}$. Hence if $A$ and $\tilde A$ are isomorphic, we obtain that $S_{\pi}$, $S_{\tilde\pi}$ are affinely equivalent as required.

Conversely, let ${\mathcal A}:{\mathfrak m}\ra\tilde{\mathfrak m}$ be an affine equivalence with ${\mathcal A}(S_{\pi})=S_{\tilde\pi}$, and $y:={\mathcal A}(0)$. Consider the linear map ${\mathcal L}(u):={\mathcal A}(u)-y$, with $u\in{\mathfrak m}$. We will show that ${\mathcal L}:{\mathfrak m}\ra\tilde{\mathfrak m}$ is an algebra isomorphism, which will imply that $A$ and $\tilde A$ are isomorphic. 

Clearly, ${\mathcal L}$ maps $S_{\pi}$ onto $S_{\tilde\pi}-y$. Consider the admissible projection on $\tilde A$ given by the formula $\tilde\pi'(a):=\tilde\pi(\widetilde{\exp}(y)a)$, where $\widetilde{\exp}$ is the exponential map associated to $\tilde A$. We then have $S_{\tilde\pi}-y=S_{\tilde\pi'}$, hence ${\mathcal L}$ maps $S_{\pi}$ onto $S_{\tilde\pi'}$. 

Recall that $S_{\pi}$ is the graph of the polynomial map $P_{\pi}:{\mathcal K}\ra\Ann({\mathfrak m})$ (see (\ref{poly})). Set $n:=\dim{\mathfrak m}-1=\dim\tilde{\mathfrak m}-1$ and choose coordinates\linebreak $\alpha=(\alpha_1,\dots,\alpha_n)$ in ${\mathcal K}$ and a coordinate $\alpha_{n+1}$ in $\Ann({\mathfrak m})$. In these coordinates the hypersurface $S_{\pi}$ is written as
$$
\alpha_{n+1}=\sum_{i,j=1}^ng_{ij}\alpha_i\alpha_j+\sum_{i,j,\ell=1}^nh_{ijk}\alpha_i\alpha_j\alpha_{\ell}+\dots,
$$ 
where $g_{ij}$ and $h_{ij\ell}$ are symmetric in all indices, and the dots denote the higher-order terms. Furthermore, $(g_{ij})$ is non-degenerate since the form $b_{\pi}$ defined in (\ref{bpi}) is non-degenerate.  In Proposition 2.10 in \cite{FIKK} (which works over any field of zero characteristic) we showed that the above equation of $S_{\pi}$ is in {\it Blaschke normal form}, that is, one has $\sum_{ij=1}^ng^{ij}h_{ij\ell}=0$ for all $\ell$, where $(g^{ij}):=(g_{ij})^{-1}$.

We will now need the following lemma, which is a variant of the second statement of Proposition 1 in \cite{EE}.

\begin{lemma}\label{Blashke} \sl Let $V$, $W$ be vector spaces over an infinite field $\FF$, with\linebreak $\dim V=\dim W=N+1$, $N\ge 0$. Choose coordinates $\beta=(\beta_1,\dots,\beta_N)$, $\beta_{N+1}$ in $V$ and coordinates $\gamma=(\gamma_1,\dots,\gamma_N)$, $\gamma_{N+1}$ in $W$. Let $S\subset V$, $T\subset W$ be hypersurfaces given, respectively, by the equations
\begin{equation}
\beta_{N+1}=P(\beta_1,\dots,\beta_N),\quad \gamma_{N+1}=Q(\gamma_1,\dots,\gamma_N),\label{sim}
\end{equation}
where $P$, $Q$ are polynomials without constant and linear terms. Assume further that equations {\rm (\ref{sim})} are in Blaschke normal form. Then every bijective linear transformation of $V$ onto $W$ that maps $S$ into $T$ has the form
$$
\gamma=C\beta,\quad \gamma_{N+1}= c\,\beta_{N+1},
$$
where $C\in\GL(N,\FF)$, $c\in\FF^*$.
\end{lemma}

\noindent{\bf Proof:} Let $L:V\ra W$ be a bijective linear transformation. Write $L$ in the most general form 
$$
\gamma=C\beta+d\beta_{N+1},\quad \gamma_{N+1}=\sum_{i=1}^N c_i\beta_i+ c\,\beta_{N+1}
$$
for some $c_1,\dots,c_N,c\in\FF$, $d\in\FF^N$, and $N\times N$-matrix $C$ with entries in $\FF$. Then the condition $L(S)\subset T$ is expressed as 
\begin{equation}
\sum_{i=1}^N c_i\beta_i+cP(\beta)\equiv Q(C\beta+dP(\beta)).\label{equal4}
\end{equation} 

Since $\FF$ is an infinite field, identity (\ref{equal4}) implies that the coefficients at the same monomials on the left and on the right are equal (see, e.g. Section 1.1 in \cite{W}). Then, since $P$ and $Q$ do not contain linear terms, it follows that $c_1=\dots=c_N=0$ and therefore $C\in\GL(N,\FF)$, $c\in\FF^*$. Further, comparing the second- and third-order terms in (\ref{equal4}), it is straightforward to see that the equations of both $S$ and $T$ cannot be in Blaschke normal form unless $d=0$. \qed
\vspace{-0.1cm}\\

By Lemma \ref{Blashke}, writing the hypersurface $S_{\tilde\pi'}$ in coordinates $\tilde\alpha=(\tilde\alpha_1,\dots,\tilde\alpha_n)$, $\tilde\alpha_{n+1}$ in $\tilde{\mathfrak m}$ chosen as above, we see that the map  ${\mathcal L}$ has the form
\begin{equation}
\quad \tilde\alpha= C\alpha,\quad \tilde\alpha_{n+1}=c\,\alpha_{n+1}\label{coord}
\end{equation}
for some $C\in\GL(n,\FF)$, $c\in\FF^*$. In coordinate-free formulation, (\ref{coord})\linebreak means that with respect to the decompositions ${\mathfrak m}={\mathcal K}\oplus\Ann({\mathfrak m})$ and\linebreak $\tilde{\mathfrak m}=\tilde{\mathcal K}\oplus\Ann(\tilde{\mathfrak m})$, where $\tilde{\mathcal K}:=\ker\tilde\pi'\cap\tilde{\mathfrak m}$, the map ${\mathcal L}$ has the    block-diagonal form, that is, there exist linear isomorphisms $L_1: {\mathcal K}\ra\tilde{\mathcal K}$ and $L_2:\Ann({\mathfrak m})\ra\Ann(\tilde{\mathfrak m})$ such that ${\mathcal L}(u+a)=L_1(u)+L_2(a)$, with $u\in{\mathcal K}$, $a\in\Ann({\mathfrak m})$. Therefore, for the corresponding polynomial maps $P_{\pi}$ and $P_{\tilde\pi'}$ (see (\ref{poly})) we have
\begin{equation}
L_2\circ P_{\pi}=P_{\tilde\pi'}\circ L_1.\label{lineq1}
\end{equation}
Clearly, identity (\ref{lineq1}) yields
\begin{equation}
L_2\circ P_{\pi}^{[m]}=P_{\tilde\pi'}^{[m]}\circ L_1,\quad m=2,\dots,\nu,\label{lineq}
\end{equation}
where $P_{\pi}^{[m]}$, $P_{\tilde\pi'}^{[m]}$ are the homogeneous components of degree $m$ of $P_{\pi}$, $P_{\tilde\pi'}$, respectively, and $\nu$ is the nil-index of each of ${\mathfrak m}$ and $\tilde{\mathfrak m}$ (observe that the nil-indices of ${\mathfrak m}$ and $\tilde{\mathfrak m}$ are equal since by (\ref{lineq1}) one has $\deg P_{\pi}=\deg P_{\tilde\pi'}$). 

Let $\omega_m$ be the symmetric $\Ann({\mathfrak m})$-valued $m$-form on ${\mathcal K}$ defined as follows:
\begin{equation}
\omega_m(u_1,\dots,u_m):=\pi(u_1\dots u_m),\quad u_1,\dots,u_m\in{\mathcal K},\,\,m=2,\dots,\nu.\label{defforms}
\end{equation} 
By (\ref{poly}) we have 
\begin{equation}
P_{\pi}^{[m]}(u)=-\frac{1}{m!}\omega_m(u,\dots,u),\quad u\in{\mathcal K},\,\, m=2,\dots,\nu.\label{connect1}
\end{equation}
We will concentrate on the forms $\omega_2$ and $\omega_3$ (in fact, it is shown in Proposition 2.8 in \cite{FIKK} that every $\omega_m$ with $m>2$ is completely determined by $\omega_2$, $\omega_3$). As in \cite{FIKK}, we define a commutative product $(u,v)\mapsto u* v$ on ${\mathcal K}$ by requiring the identity
\begin{equation}
\omega_2(u* v,w)=\omega_3(u,v,w)\label{defprod}
\end{equation}
to hold for all $u,v,w\in{\mathcal K}$. Owing to the non-degeneracy of the form $b_{\pi}$ defined in (\ref{bpi}), the form $\omega_2$ is non-degenerate and therefore for any $u,v\in{\mathcal K}$ the element $u*v\in{\mathcal K}$ is uniquely determined by (\ref{defprod}). 

We need the following lemma.

\begin{lemma}\label{newprod}\sl For any two elements $u+a$, $v+b$ of ${\mathfrak m}$, with $u,v\in{\mathcal K}$ and $a,b\in\Ann({\mathfrak m})$, one has
\begin{equation}
(u+a)(v+b)=u* v+\omega_2(u,v).\label{sameproduct}
\end{equation}
\end{lemma}

\noindent {\bf Proof:} We write the product $uv$ with respect to the decomposition\linebreak ${\mathfrak m}={\mathcal K}\oplus\Ann({\mathfrak m})$ as $uv=(uv)_1+(uv)_2$, where $(uv)_1\in{\mathcal K}$ and $(uv)_2\in\Ann({\mathfrak m})$. It is then clear that $(uv)_2=\omega_2(u,v)$. Therefore, to prove (\ref{sameproduct}) we need to show that $u*v=(uv)_1$. This identity is proved by a simple calculation similar to the one that occurs in the proof of Proposition 2.8 of \cite{FIKK} (cf. Proposition 5.13 in \cite{FK}). Indeed, from (\ref{defforms}), (\ref{defprod}) for any $w\in{\mathcal K}$ one has
$$
\begin{array}{l}
\omega_2(u*v,w)=\omega_3(u,v,w)=\pi(uvw)=\\
\vspace{-0.1cm}\\
\hspace{4cm}\pi\Bigl([(uv)_1+(uv)_2]w\Bigr)=\pi((uv)_1w)=\omega_2((uv)_1,w).
\end{array}
$$
Since $\omega_2$ is non-degenerate, the above identity implies $u*v=(uv)_1$ as required.\qed
\vspace{0.1cm}\\

Let $\tilde\omega_m$ be the symmetric $\Ann(\tilde{\mathfrak m})$-valued forms on $\tilde{\mathcal K}$ arising from the admissible projection $\tilde\pi'$ on $\tilde{\mathfrak m}$ as in (\ref{defforms}), with $m=2,\dots,\nu$. By (\ref{lineq}), (\ref{connect1}) we have
\begin{equation}
\begin{array}{l}
L_2(\omega_m(u_1,\dots,u_m))=\tilde\omega_m(L_1(u_1),\dots, L_1(u_m)),\\ 
\vspace{-0.1cm}\\
\hspace{6cm}u_1,\dots,u_m\in{\mathcal K},\,\,m=2,\dots,\nu.
\end{array}\label{impidentities}
\end{equation}
Denote by $\tilde *$ the product on $\tilde{\mathcal K}$ defined by $\tilde\omega_2$, $\tilde\omega_3$ as in (\ref{defprod}). Now, for all $u,v\in{\mathcal K}$ and $a,b\in\Ann({\mathfrak m})$, Lemma \ref{newprod} and identity (\ref{impidentities}) with $m=2$ yield
\begin{equation}
\makebox[250pt]{$\begin{array}{l}
{\mathcal L}\Bigl((u+a)(v+b)\Bigr)={\mathcal L}(u*v+\omega_2(u,v))=L_1(u*v)+L_2(\omega_2(u,v)),\\
\vspace{0cm}\\
{\mathcal L}(u+a){\mathcal L}(v+b)=(L_1(u)+L_2(a))(L_1(v)+L_2(b))=\\
\vspace{-0.3cm}\\
\hspace{2cm}L_1(u)\tilde *L_1(v)+\tilde\omega_2(L_1(u),L_1(v))=L_1(u)\tilde *L_1(v)+L_2(\omega_2(u,v)).
\end{array}$}\label{imporeq}
\end{equation}
Next, for any $w\in{\mathcal K}$, from (\ref{defprod}) and identity (\ref{impidentities}) with $m=2,3$ one obtains
$$
\begin{array}{l}
\tilde\omega_2(L_1(u)\tilde*L_1(v),L_1(w))=\tilde\omega_3(L_1(u),L_1(v),L_1(w))=\\
\vspace{-0.1cm}\\
\hspace{3cm}L_2(\omega_3(u,v,w))=L_2(\omega_2(u*v,w))=\tilde\omega_2(L_1(u*v),L_1(w)),
\end{array}
$$
which implies $L_1(u)\tilde*L_1(v)=L_1(u*v)$. It then follows from (\ref{imporeq}) that ${\mathcal L}\Bigl((u+a)(v+b)\Bigr)={\mathcal L}(u+a){\mathcal L}(v+b)$, that is, ${\mathcal L}:{\mathfrak m}\ra\tilde{\mathfrak m}$ is an algebra isomorphism. The proof of the theorem is complete.\qed

\section{Example of application of Theorem \ref{equivalence}}\label{app}
\setcounter{equation}{0}

Theorem \ref{equivalence} is particularly useful when at least one of the hypersurfaces $S_{\pi}$, $S_{\tilde\pi}$ is affinely homogeneous (recall that a subset ${\mathcal S}$ of a vector space $V$ is called affinely homogeneous if for every pair of points $p,q\in{\mathcal S}$ there exists a bijective affine map ${\mathcal A}$ of $V$ such that ${\mathcal A}({\mathcal S})={\mathcal S}$ and ${\mathcal A}(p)=q$). In this case the hypersurfaces $S_{\pi}$, $S_{\tilde\pi}$ are affinely equivalent if and only if they are linearly equivalent. Indeed, if, for instance, $S_{\pi}$ is affinely homogeneous and ${\mathcal A}:{\mathfrak m}\ra\tilde{\mathfrak m}$ is an affine equivalence between $S_{\pi}$, $S_{\tilde\pi}$, then ${\mathcal A}\circ {\mathcal A}'$ is a linear equivalence between $S_{\pi}$, $S_{\tilde\pi}$, where ${\mathcal A}'$ is an affine automorphism of $S_{\pi}$ such that ${\mathcal A}'(0)={\mathcal A}^{-1}(0)$. Clearly, in this case $S_{\tilde\pi}$ is affinely homogeneous as well. 

Now, as we observed in the proof of Theorem \ref{equivalence}, every linear equivalence ${\mathcal L}$ between $S_{\pi}$, $S_{\tilde\pi}$ has the block-diagonal form with respect to the decompositions ${\mathfrak m}={\mathcal K}\oplus\Ann({\mathfrak m})$ and $\tilde{\mathfrak m}=\tilde{\mathcal K}\oplus\Ann(\tilde{\mathfrak m})$, where $\tilde{\mathcal K}:=\ker\tilde\pi\cap\tilde{\mathfrak m}$, that is, there exist linear isomorphisms $L_1: {\mathcal K}\ra\tilde{\mathcal K}$ and $L_2:\Ann({\mathfrak m})\ra\Ann(\tilde{\mathfrak m})$ such that ${\mathcal L}(u+a)=L_1(u)+L_2(a)$, with $u\in{\mathcal K}$, $a\in\Ann({\mathfrak m})$. Therefore, analogously to (\ref{lineq}), for the corresponding polynomial maps $P_{\pi}$ and $P_{\tilde\pi}$ (see (\ref{poly})) we have
\begin{equation}
L_2\circ P_{\pi}^{[m]}=P_{\tilde\pi}^{[m]}\circ L_1\quad\hbox{for all $m\ge 2$},\label{lineq2}
\end{equation}
where, as before, $P_{\pi}^{[m]}$, $P_{\tilde\pi}^{[m]}$ are the homogeneous components of degree $m$ of $P_{\pi}$, $P_{\tilde\pi}$, respectively.

Thus, Theorem \ref{equivalence} yields the following corollary (cf. Theorem 2.11 in \cite{FIKK}) .

\begin{corollary}\label{cor}\sl Let $A$, $\tilde A$ be Gorenstein algebras of finite vector space dimension greater than 1, and $\pi$, $\tilde\pi$ admissible projections on $A$, $\tilde A$, respectively.

\noindent {\rm (i)} If $A$ and $\tilde A$ are isomorphic and at least one of $S_{\pi}$, $S_{\tilde\pi}$ is affinely homogeneous, then for some linear isomorphisms $L_1: {\mathcal K}\ra\tilde{\mathcal K}$ and $L_2:\Ann({\mathfrak m})\ra\Ann(\tilde{\mathfrak m})$ identity {\rm (\ref{lineq2})} holds. In this case both $S_{\pi}$ and $S_{\tilde\pi}$ are affinely homogeneous.

\noindent {\rm (ii)} If for some linear isomorphisms $L_1: {\mathcal K}\ra\tilde{\mathcal K}$ and $L_2:\Ann({\mathfrak m})\ra\Ann(\tilde{\mathfrak m})$ identity {\rm (\ref{lineq2})} holds, then the hypersurfaces $S_{\pi}$, $S_{\tilde\pi}$ are linearly equivalent and therefore the algebras $A$ and $\tilde A$ are isomorphic. 
\end{corollary}

In \cite{I2} (see also \cite{FK}) we found a criterion for the affine homogeneity of some (hence every) hypersurface $S_{\pi}$ arising from a Gorenstein algebra $A$. Namely, $S_{\pi}$ is affinely homogeneous if and only if the action of the automorphism group of the algebra ${\mathfrak m}$ on the set of all hyperplanes in ${\mathfrak m}$ complementary to $\Ann({\mathfrak m})$ is transitive. Furthermore, we showed that this condition is satisfied if $A$ is non-negatively graded in the sense that it can be represented as a direct sum 
\begin{equation}
A=\bigoplus_{j\ge0}A^j,\quad A^{j}A^{\ell}\subset A^{j+\ell},\label{sum}
\end{equation}
where $A^{j}$ are linear subspaces of $A$, with $A^0\simeq\FF$ (in this case ${\mathfrak m}=\oplus_{j>0}A^j$ and $\Ann({\mathfrak m})=A^d$ for $d:=\max\{j:A^{j}\ne 0\}$). It then follows that part (i) of Corollary \ref{cor} applies in the situation when one (hence the other) of the algebras $A$, $\tilde A$ is non-negatively graded. Note, however, that the existence of a non-negative grading on $A$ is not a necessary condition for the affine homogeneity of $S_{\pi}$ (see, e.g. Remark 2.6 in \cite{FIKK}). Also, as shown in Section 8.2 in \cite{FK}, the hypersurface $S_{\pi}$ need not be affinely homogeneous in general.

To demonstrate how our method works, we will now apply Corollary \ref{cor} to a one-parameter family of non-negatively graded Gorenstein algebras. For $t\in\FF$, $t\ne\pm 2$, define
$$
A_t:=\FF[x,y]/(2x^3+txy^3,tx^2y^2+2y^5).
$$
It is straightforward to verify that every $A_t$ is a Gorenstein algebra of dimension 15. We will prove the following proposition.

\begin{proposition}\label{solution}\sl $A_r$ and $A_s$ are isomorphic if and only if $r=\pm s$.
\end{proposition}

\noindent {\bf Proof:} The sufficiency implication is trivial (just replace $y$ by $-y$). For the converse implication, consider the following monomials in $\FF[x,y]$:
$$
1,\,\,x,\,\,y,\,\,x^2,\,\,xy,\,\,y^2,\,\,x^2y,\,\,xy^2,\,\,y^3,\,\,xy^3,\,\,x^2y^2,\,\,y^4,\,\,x^2y^3,\,\,xy^4,\,\,x^2y^4.
$$
Let $e_0,\dots,e_{14}$, respectively, be the vectors in $A_t$ arising from these monomials (to simplify the notation, we do not indicate the dependence of $e_j$ on $t$). They form a basis of $A_t$. Define
$$
\begin{array}{l}
A^0_t:=\langle e_0\rangle,\, A^1_t:=0,\, A^2_t:=\langle e_2\rangle,\, A^3_t:=\langle e_1\rangle,\, A^4_t:=\langle e_5\rangle,\\
\vspace{-0.1cm}\\
A^5_t:=\langle e_4\rangle,\,A^6_t:=\langle e_3,e_8\rangle,\,A^7_t:=\langle e_7\rangle,\,A^8_t:=\langle e_6, e_{11}\rangle,\,A^9_t:=\langle e_9\rangle,\\
\vspace{-0.1cm}\\
A^{10}_t:=\langle e_{10}\rangle,\,A^{11}_t:=\langle e_{13}\rangle,\,A^{12}_t:=\langle e_{12}\rangle,\,A^{13}_t:=0,\,A^{14}_t:=\langle e_{14}\rangle,
\end{array}
$$ 
where $\langle{}\cdot{}\rangle$ denotes linear span. It is straightforward to check that the subspaces $A^j_t$ form a non-negative grading on $A_t$.  

Next, denote by ${\mathfrak m}_t$ the maximal ideal of $A_t$ and let $\pi_t$ be the projection on $A_t$ with range $\Ann({\mathfrak m}_t)=A^{14}_t$ such that $\ker\pi_t=\langle e_0,\dots,e_{13}\rangle$. Denote by $\alpha_1,\dots,\alpha_{14}$ the coordinates in ${\mathfrak m}_t$ with respect to the basis $e_1,\dots,e_{14}$. In these coordinates the corresponding polynomial map $P_t:=P_{\pi_{{}_t}}$ is written as  
$$
\begin{array}{l}
\displaystyle P_t=\frac{t}{10080}\alpha_2^7-\frac{1}{48}\alpha_2^4\left(\alpha_1^2-\frac{t}{5}\alpha_2\alpha_5\right)+\frac{t}{48}\alpha_1^4\alpha_2-\frac{1}{4}\alpha_1^2\alpha_2^2\alpha_5-\frac{1}{6}\alpha_1\alpha_2^3\alpha_4+\\
\vspace{-0.1cm}\\
\displaystyle\hspace{1cm}\frac{t}{24}\alpha_2^3\alpha_5^2+\frac{t}{48}\alpha_2^4\alpha_8-\frac{1}{24}\alpha_2^4\alpha_3+\hbox{terms of degree $\le 4\,$}.
\end{array}
$$

Suppose that for some $r\ne s$ the algebras $A_r$ and $A_s$ are isomorphic. By part (i) of Corollary \ref{cor} there exist $C\in\GL(13,\FF)$ and $c\in\FF^*$ such that   
\begin{equation}
c\, P_r(\alpha)\equiv P_s(C\alpha)\,.\label{equivnew}
\end{equation}
Since 0 is the only value of $t$ for which $P_t$ has degree 7, we have $r,s\ne 0$. Comparing the terms of order 7 in identity (\ref{equivnew}), we obtain that the second row in the matrix $C$ has the form $(0,\mu,0,\dots,0)$, and
\begin{equation}
c=\frac{s}{r}\mu^7\,.\label{7}
\end{equation}
Next, comparing the terms of order 6 in (\ref{equivnew}), we see that the first row in the matrix $C$ has the form $(\sigma,\rho,0,\dots,0)$,  and 
\begin{equation}
c=\mu^4\sigma^2\,.\label{6}
\end{equation}
Further, comparing the terms of order 5 in (\ref{equivnew}) that do not involve $\alpha_2^2$, we obtain
\begin{equation}
c=\frac{s}{r}\mu\sigma^4\,.\label{5}
\end{equation}
Now, (\ref{7}), (\ref{6}), (\ref{5}) yield $r^2=s^2$ as required. \qed
\vspace{0.3cm}

One can give an alternative proof of Proposition \ref{solution} by utilizing Maca\-ulay's inverse systems. For example, the polynomial
$$
\mu_t:=\frac{t}{10080}y^7-\frac{1}{48}x^2y^4+\frac{t}{48}x^4y
$$
is an inverse system for $A_t$. Although $\mu_t$ are simpler than $P_t$, the equivalence relation that one would have to use for them (see formula (7) in \cite{ER}) is much more complicated than linear equivalence up to scale, which we utilized for $P_t$. Note also that upon identification of $\alpha_1,\alpha_2$ with $x,y$, respectively, $\mu_t$ is obtained from $P_t$ by setting $\alpha_3=\dots=\alpha_{13}=0$. In the next section we will see that this fact is a manifestation of a general principle.

\section{Nil-polynomials and Macaulay's\\ inverse systems}\label{invsys}
\setcounter{equation}{0}

As before, let $A$ be an Artinian Gorenstein algebra with $\dim_{\FF}A>1$ and maximal ideal ${\mathfrak m}$ of nil-index $\nu$. Fix a hyperplane $\Pi$ in ${\mathfrak m}$ complementary to $\Ann({\mathfrak m})$. An $\FF$-valued polynomial $P$ on $\Pi$ is called a {\it nil-polynomial}\, if there exists a linear form $\omega:A\ra\FF$ such that $\ker\omega=\langle\Pi,{\bf 1}\rangle$ and $P=\omega\circ\exp|_{\Pi}$. Observe that $\deg P=\nu$. Upon identification of $\Ann({\mathfrak m})$ with $\FF$, the class of nil-polynomials coincides with that of $\Ann({\mathfrak m})$-valued polynomial maps $P_{\pi}$ introduced in (\ref{poly}). 

Next, let $k:=\emb\dim A:=\dim_{\FF}{\mathfrak m}/{\mathfrak m}^2\ge 1$ be the embedding dimension of $A$. We now assume that $\dim_{\FF}A>2$. Then $\nu\ge 2$ and $\Pi$ contains a $k$-dimensional subspace that forms a complement to ${\mathfrak m}^2$ in ${\mathfrak m}$. Fix any such subspace $L$, choose a basis $e_1,\dots,e_k$ in $L$ and let $y_1,\dots,y_k$ be the coordinates in $L$ with respect to this basis. Denote by $Q\in\FF[y_1,\dots,y_k]$ the restriction of the nil-polynomial $P$ to $L$ written in these coordinates. Clearly, one has
$$
Q(y_1,\dots,y_k)=\sum_{j=0}^{\nu}\frac{1}{j!}\omega\Bigl((y_1e_1+\dots+y_ke_k)^j\Bigr).
$$

Further, the elements $e_1,\dots,e_k$ generate $A$ as an algebra, hence $A$ is isomorphic to $\FF[x_1,\dots,x_k]/I$, where $I$ is the ideal of all relations, i.e. polynomials $f\in\FF[x_1,\dots,x_k]$ with $f(e_1,\dots,e_k)=0$. Observe that $I$ contains the monomials $x_1^{\nu+1},\dots,x_k^{\nu+1}$. It is well-known that, since the quotient $\FF[x_1,\dots,x_k]/I$ is Gorenstein, there is a polynomial $g\in\FF[y_1,\dots,y_k]$ of degree $\nu$ such that $I=\Ann(g)$, where
$$
\Ann(g):=\left\{f\in\FF[x_1,\dots,x_k]: f\left(\frac{\partial}{\partial y_1},\dots,\frac{\partial}{\partial y_k}\right)(g)=0\right\}
$$
is the annihilator of $g$ (see, e.g. \cite{ER} and references therein). Any such polynomial $g$ is called {\it Macaulay's inverse system}\, for the Gorenstein quotient $\FF[x_1,\dots,x_k]/I$. It is well-known that inverse systems can be used for solving the isomorphism problem for quotients of this kind (see, e.g. Proposition 2.2 in \cite{ER}).  

We are now ready to state the result of this section.

\begin{theorem}\label{main}\sl The polynomial $Q$ is an inverse system for the quotient $\FF[x_1,\dots,x_k]/I$.
\end{theorem}

\noindent{\bf Proof:} Fix any polynomial $f\in\FF[x_1,\dots,x_k]$
$$
f=\sum_{0\le i_1,\dots,i_k\le N}a_{i_1,\dots,i_k}x_1^{i_1}\dots x_k^{i_k}
$$ 
and calculate
$$
\displaystyle f\left(\frac{\partial}{\partial y_1},\dots,\frac{\partial}{\partial y_k}\right)(Q)=\sum_{0\le i_1,\dots,i_k\le N}a_{i_1,\dots,i_k}\sum_{j= i_1+\dots+i_k}^{\nu}\frac{1}{(j-(i_1+\dots+i_k))!}\times
$$
\begin{equation}
\begin{array}{l}
\displaystyle\omega\Bigl((y_1e_1+\dots+y_ke_k)^{j-(i_1+\dots+i_k)}e_1^{i_1}\dots e_k^{i_k}\Bigr)=\\
\vspace{-0.3cm}\\
\displaystyle\sum_{m=0}^{\nu}\frac{1}{m!}\omega\Bigl((y_1e_1+\dots+y_ke_k)^m\sum_{\tiny\begin{array}{l} 0\le i_1,\dots,i_k\le N,\\\vspace{-0.1cm}\\ i_1+\dots+i_k\le\nu-m\end{array}}a_{i_1,\dots,i_k}e_1^{i_1}\dots e_k^{i_k}\Bigr)=\\
\vspace{-0.3cm}\\
\displaystyle\sum_{m=0}^{\nu}\frac{1}{m!}\omega\Bigl((y_1e_1+\dots+y_ke_k)^m\,f(e_1,\dots,e_k)\Bigr).
\end{array}\label{diff}
\end{equation}
Formula (\ref{diff}) immediately implies that $I\subset\Ann(Q)$.

Conversely, let $f\in\FF[x_1,\dots,x_k]$ be an element of $\Ann(Q)$. Then (\ref{diff}) yields
\begin{equation}
\sum_{m=0}^{\nu}\frac{1}{m!}\omega\Bigl((y_1e_1+\dots+y_ke_k)^m\,f(e_1,\dots,e_k)\Bigr)=0.\label{eq2}
\end{equation}
Collecting the terms containing $y_1^{i_1}\dots y_k^{i_k}$ in (\ref{eq2}) we obtain
\begin{equation}
\omega\Bigl(e_1^{i_1}\dots e_k^{i_k}\,f(e_1,\dots,e_k)\Bigr)=0\label{eq3}
\end{equation}
for all indices $i_1,\dots,i_k$. Since $e_1,\dots,e_k$ generate $A$, identities (\ref{eq3}) yield 
\begin{equation}
\omega\Bigl(A\, f(e_1,\dots,e_k)\Bigr)=0.\label{nondeg}
\end{equation}
Further, since the bilinear form $(a,b)\mapsto \omega(ab)$ is non-degenerate on $A$, identity (\ref{nondeg}) implies $f(e_1,\dots,e_k)=0$. Therefore $f\in I$, which shows that $I=\Ann(Q)$ as required.\qed

\begin{remark}\label{stgradedhom}\rm Suppose that $A$ is a standard graded algebra, i.e. $A$ can be represented in the form (\ref{sum}) with $A^j=(A^1)^j$. Choose
$$
\Pi:=\bigoplus_{j=1}^{\nu-1}A^j,\quad L:=A^1.
$$
In this case the ideal $I$ is homogeneous, i.e. generated by homogeneous relations. For any nil-polynomial $P$ on $\Pi$ its restriction $Q$ to $L$ coincides with the homogeneous component of degree $\nu$ of $P$:
$$
Q(y_1,\dots,y_k)=\frac{1}{\nu!}\omega\Bigl((y_1e_1+\dots+y_ke_k)^{\nu}\Bigr),
$$
which is a well-known way of writing inverse systems for Gorenstein quotients by homogeneous ideals.

\end{remark}

{\obeylines
\noindent Department of Mathematics
\noindent The Australian National University
\noindent Canberra, ACT 0200
\noindent Australia
\noindent e-mail: alexander.isaev@anu.edu.au
}

\end{document}